\documentclass[a4 paper,12pt]{amsart}
\usepackage{amsmath,amsthm,amssymb,latexsym,graphicx}
\usepackage{hyperref}
\usepackage{xypic}

\newtheorem{theorem}{Theorem} [section]

\newtheorem{proposition}[theorem]{Proposition}

\theoremstyle{definition}

\newtheorem*{1'}{Definition 1$'$}

\theoremstyle{remark}
\newtheorem*{remark}{Remark}

\begin{document}
\title{On The Cofibrant Generation Of Model Categories}
\author{George Raptis}
\address{Mathematical Institute\\
	24-29 St Giles'\\
	Oxford\\
	OX1 3LB\\
	England}
\email{raptis@maths.ox.ac.uk}
\date{}
\maketitle

\begin{abstract}
The paper studies the problem of the cofibrant generation of a model category.  We prove that, assuming Vop\v{e}nka's principle, every cofibrantly generated model category is Quillen equivalent to a combinatorial model category. We discuss cases where this result implies that the class of 
weak equivalences in a cofibrantly generated model category is accessibly embedded. We also prove a necessary condition for a model category to be cofibrantly generated by a set of generating cofibrations between cofibrant objects.  
\end{abstract}

\section{Introduction and Statement of Results}

\parindent=0.2in The purpose of the paper is to study the problem of the cofibrant generation of a model category and relate it to ideas from the theory of combinatorial model categories.  A combinatorial model category is a cofibrantly generated model category whose underlying category is locally presentable.  They were first introduced by J. Smith and have been studied fruitfully ever since \cite{Be} \cite{D1,D2}, \cite{Ro}. Locally presentable categories allow empowered uses of Quillen's  small-object argument that greatly facilitate the construction of combinatorial model structures. In fact, the problem of the existence of a model category structure on a locally presentable category often reduces to the problem of the cofibrant generation for the candidate class of cofibrations, as long as the class of weak equivalences is known to satisfy some closure and smallness conditions. Moreover, combinatorial model categories share remarkable categorical and homotopical properties and include many important examples. Simplicial sets and, more generally Grothendieck topoi, are locally presentable categories that carry combinatorial model structures \cite{Ci}.
\

\parindent=0.2in Whereas not every model category is combinatorial, most of the important examples are at least Quillen equivalent to one. Assuming some input from set theory, one of our main results formalises this claim as follows.

\begin{theorem} \label{A}
Assuming Vop\v{e}nka's principle, every cofibrantly generated model category is Quillen equivalent to a combinatorial model category.
\end{theorem}

\parindent=0.2in Vop\v{e}nka's principle is a set-theoretical axiom that characteristically appears in the study of locally presentable categories because it gives a simple characterisation of them. Under the assumption that Vop\v{e}nka's principle holds, a category is locally presentable if and only if it is cocomplete and has a dense subcategory [\cite{AR}, Theorem 6.14]. We will discuss cases where the assumption of Vop\v{e}nka's principle in the theorem is not needed.  

\parindent=0.2in In the most general case, the problem of the cofibrant generation for a cofibrantly closed class of morphisms $\mathcal{S}$ in a cocomplete category $\mathcal{C}$ asks whether there exists a set of morphisms $\mathcal{I}$ whose cofibrant closure in $\mathcal{C}$ is $\mathcal{S}$. This is difficult to decide in this generality.  The problem is especially interesting in the following two closely related cases: 
\begin{itemize}
\item[$\mathbf{A}$] When is a model category $\mathcal{M}$ cofibrantly generated?
\item[$\mathbf{B}$] Let $\mathcal{C}$ be a locally presentable category with a class of weak equivalences $\mathcal{W}$ which satisfies the 2-out-of-3 property, it is closed under retracts and it is accessible and accessibly embedded in $\mathcal{C}^{\rightarrow}$.  Let $\mathcal{C}of$ be a cofibrantly closed class of maps in $\mathcal{C}$ such that $\mathcal{C}of \cap \mathcal{W}$ is closed under pushouts and transfinite compositions. When is $\mathcal{C}of$ cofibrantly generated? 
\end{itemize}

\parindent=0.2in   Regarding the first question, $\mathcal{M}$ is cofibrantly generated only if there is a set of "test"-cofibrations (resp. trivial "test"-cofibrations) such that every morphism is a trivial fibration (resp. fibration) if and only if it has the right lifting property with respect to this set. Cofibrantly generated model categories include most examples of interest in applications and they allow certain constructions to be possible directly (e.g. model structures on diagram categories, constructions of homotopy colimits, etc. see \cite{Hi}) by essentially giving a grip on the cofibrations the same way that CW complexes do in the homotopy theory of topological spaces. 

\parindent=0.2in A set $S$ of objects in a category $\mathcal{C}$ is called left adequate if, for every map $f: X \to Y$, $f$ is an isomorphism if and only if $\mathcal{C}(K,X) \to \mathcal{C}(K,Y)$ is a bijection for all $K \in S$. This concept is due to Heller \cite{He}. A set $S$ of objects in a category $\mathcal{C}$ with a terminal object is called left weakly adequate if, for every object $X$ in $\mathcal{C}$, $X$ is isomorphic to the terminal object if and only if $\mathcal{C}(K,X) = \{ \star \}$ for every $K \in S$.\

\parindent=0.2in We prove the following necessary smallness condition. This was previously known for cofibrantly generated $pointed$ model categories [\cite{Ho}, Theorem 7.3.1].

\begin{theorem} \label{B}
Let $\mathcal{M}$ be a cofibrantly generated model category. Suppose that there is a set $I$ of generating cofibrations between cofibrant objects. Then the homotopy category $Ho \mathcal{M}$ of $\mathcal{M}$ admits a left weakly adequate set of objects.
\end{theorem}

\parindent=0.2in  Not every model category is cofibrantly generated and examples are known in the literature, e.g. \cite{Ch}. In \cite{St}, Str\o{}m discovered a model structure on the category of topological spaces whose weak equivalences are the homotopy equivalences. We conjecture that Theorem \ref{B} applies to show that Str\o{}m's model category is not cofibrantly generated.

\parindent=0.2in Regarding the second question above, it is, in practice, a crucial step in order to deduce that $\mathcal{C}of$ and  $\mathcal{W}$ determine a combinatorial model category structure on $\mathcal{C}$. This observation rests on J. Smith's main theorem, a version of which will be recalled below.

\parindent=0.2in The organisation of the paper is as follows. In Section 2, we recall the definitions of cofibrantly generated and combinatorial model categories. Its purpose is mainly to establish some notation and terminology. In this we follow the conventions of \cite{Hi},\cite{Ho} that we recommend.  For background in the theory of locally presentable categories and accessibility, the reader should consult \cite{AR}. A nice exposition of the theory of combinatorial model categories can be found in \cite{Ro}. 
\

 \parindent=0.2in In Section 3, we prove Theorem \ref{A} and we discuss cases where the class of weak equivalences in a cofibrantly generated model category is accessibly embedded. Finally, in Section 4, we prove Theorem \ref{B}.

\section{Recollections}

\parindent=0.2in Let $\mathcal{C}$ be a cocomplete category and $I$ a set of morphisms. The cofibrant closure $Cof(I)$ of $I$ (in $\mathcal{C}$) is the smallest class of morphisms that contains $I$ and is closed under retracts, pushouts and transfinite compositions.  A class of morphisms $\mathcal{S}$ in $\mathcal{C}$ is cofibrantly generated if there exists a set $I$ of morphisms in $\mathcal{C}$  such that $\mathcal{S}=Cof(I)$.

\

\parindent=0.2in A $model$ category is a category $\mathcal{M}$ together with three classes of morphisms called weak equivalences, fibrations and cofibrations, each of which is closed under composition and contains the identities and satisfies the following axioms:
\begin{enumerate}
\item[M1] (Limits) All small limits and colimits exist in $\mathcal{M}$.
\item[M2] (2-out-of-3) If two out of three morphisms $f$,$g$ and $f \circ g$ are weak equivalences, then so is the third.
\item[M3] (Retracts) If $f$ is a retract of $g$ and $g$ is a fibration, cofibration or weak equivalence, then so is $f$.
\item[M4] (Lifting) Given a commutative diagram in $\mathcal{M}$
\begin{displaymath}
\xymatrix{
A \ar[r]^f \ar[d]^i & X \ar[d]^p \\
B \ar[r]^g &Y
}
\end{displaymath}
then there exists a lift $h: B \rightarrow X$ if either:
\begin{itemize}
\item
$i$ is a cofibration and $p$ is a trivial fibration (i.e. a fibration and a weak equivalence), or
\item
$i$ is a trivial cofibration (i.e. a cofibration and a weak equivalence) and $p$ is a fibration. 
\end{itemize}
\item[M5] (Factorisation)
Every morphism $f$ in $\mathcal{M}$ can be factorised functorially in the following two ways:
\begin{itemize}
\item
$f=qi$, where $i$ is a cofibration and $q$ is a trivial fibration
\item
$f=pj$, where $j$ is a trivial cofibration and $p$ is a fibration.
\end{itemize}
\end{enumerate}
\

\parindent=0.2in  A model category $\mathcal{M}$ is $cofibrantly$ $generated$ if there exist sets of morphisms $I$ and $J$ such that the following hold:
\begin{enumerate}
\item 
the domains of $I$ are small relative to $I-cellular$ morphisms;
the domains of $J$ are small relative to $J-cellular$ morphisms.
\item
the fibrations are the $J-injective$ morphisms;
the trivial fibrations are the $I-injective$ morphisms.
\end{enumerate}
\
\parindent=0.2in $I$ is called a set of generating cofibrations and $J$ a set of generating trivial cofibrations. Equivalently, a model category is cofibrantly generated if the classes of cofibrations and trivial cofibrations are cofibrantly generated and the required smallness condition of the definition is satisfied. Note that the smallness condition is automatically satisfied when the underlying category consists solely of presentable objects. 
\

\parindent=0.2in A model category $\mathcal{M}$ is $combinatorial$ if it is cofibrantly generated and its underlying category is locally presentable. Recall that a cocomplete category $\mathcal{C}$ is locally presentable if, for some regular cardinal $\lambda$, it has a set $S$ of $\lambda$-presentable objects such that every object is a $\lambda$-directed colimit of objects in $S$ \cite{AR}.   
\
 
\parindent=0.2in The tools for comparison between model categories are called Quillen functors; these are functors which preserve part of the structure and they have good derivability properties. Let $\mathcal{M}$ and $\mathcal{N}$ be model categories. A functor $F: \mathcal{M} \rightarrow \mathcal{N}$ is a left Quillen functor if $F$ is a left adjoint and it preserves cofibrations and trivial cofibrations. Equivalently, if its right adjoint $G: \mathcal{N} \to \mathcal{M}$ is a right Quillen functor, i.e. it preserves fibrations and trivial fibrations. The adjunction $F: \mathcal{M} \leftrightarrows \mathcal{N}:G$ is called a Quillen adjunction.
\

\parindent=0.2in It can be shown that left Quillen functors preserve weak equivalences between cofibrant objects. The subcategory of cofibrant objects is a (left) "deformation retract" of the whole category because, by the factorisation axiom, every object has a functorial cofibrant replacement. It follows that every left Quillen functor $F: \mathcal{M} \to \mathcal{N}$ admits a total left derived functor $LF: Ho \mathcal{M} \to Ho \mathcal{N}$ (and dually, every right Quillen functor $G$ admits a total right derived functor $RG$). Moreover, every Quillen adjunction $F:\mathcal{M} \leftrightarrows \mathcal{N}:G$ induces a derived adjunction $LF:Ho\mathcal{M} \leftrightarrows Ho\mathcal{N}:RG$ [\cite{Ho}, Chapter 1]. 
\

\parindent=0.2in A Quillen adjunction $F:\mathcal{M} \leftrightarrows \mathcal{N}:G$ is called a Quillen equivalence if $LF:Ho\mathcal{M} \leftrightarrows Ho\mathcal{N}: RG$ is an adjoint equivalence of categories.

\section{The Proof of Theorem \ref{A}}

\parindent=0.2in Let $S$ be a set of objects in a cocomplete category $\mathcal{C}$.  Denote by $\mathcal{S}$ the full subcategory of $\mathcal{C}$ with objects in $S$. Every object $X$ in $\mathcal{C}$ determines a comma category $\mathcal{S} \downarrow X$ and a canonical forgetful diagram $\underline{X}: \mathcal{S} \downarrow X \to \mathcal{C}$ with respect to $\mathcal{S}$.  Let $\kappa_{\mathcal{S}}(X)$ denote the colimit of this diagram in $\mathcal{C}$. There is an induced canonical morphism $\eta_X: \kappa_{\mathcal{S}}(X) \to X$ in $\mathcal{C}$ which is natural in $X$. We say that $X$ is $S$-generated if $\eta_X$ is an isomorphism, cf. \cite{FR}. Let $\mathcal{C}_{\mathcal{S}}$ denote the full subcategory of $S$-generated objects. 

\begin{proposition} \label{proposition}
 Let $\mathcal{C}$ be a cocomplete category and $S$ a set of objects. 
\begin{itemize}
 \item[(i)] The functor $\kappa_{\mathcal{S}}: \mathcal{C} \to \mathcal{C}_{\mathcal{S}}$ is right adjoint to the inclusion functor.
\item[(ii)] $\mathcal{C}_{\mathcal{S}}$ is locally presentable if and only if every object $A \in S$ is $\lambda$-presentable in $\mathcal{C}_{\mathcal{S}}$ for some regular ordinal $\lambda$. 
\item[(iii)] Assuming Vop\v{e}nka's principle, $\mathcal{C}_{\mathcal{S}}$ is locally presentable.
\end{itemize}

\end{proposition}
\begin{proof} (i) The adjunction isomorphism follows from the universal property of the morphism $\eta_X$. (ii) The ``only if'' is obvious. For the ``if'' part, note the $S$ is a strong generator of $\lambda$-presentable objects in $\mathcal{C}_{\mathcal{S}}$. (iii) $\mathcal{C}_{\mathcal{S}}$ is cocomplete and it has a dense subcategory.  
\end{proof}

\begin{proof}(of Theorem \ref{A})
Let $\mathcal{M}$ be a cofibrantly generated model category with a generating set of cofibrations $I$ and trivial cofibrations $J$. Let $S$ be the set of objects that appear as domains or codomains of $I \cup J$. Assuming Vop\v{e}nka's principle, $\mathcal{M}_{\mathcal{S}}$ is a locally presentable category. We claim that the model structure on $\mathcal{M}$ restricts to a model structure on $\mathcal{M}_{\mathcal{S}}$ and the adjunction $i: \mathcal{M}_{\mathcal{S}} \leftrightarrows \mathcal{M}: \kappa_{\mathcal{S}}$ is a Quillen equivalence. The only non-trivial part of the first claim is to show that factorisations exist in $\mathcal{M}_{\mathcal{S}}$. This is true in virtue of the fact that the factorisations given by the small-object argument can be performed in $\mathcal{M}_{\mathcal{S}}$ because $i$ preserves colimits and therefore the pushout and the directed colimit of $S$-generated objects is $S$-generated. The unit transformation $1 \to \kappa_{\mathcal{S}} \circ i$ is a natural isomorphism and therefore also a natural weak equivalence of functors. Each component of the counit $i \circ \kappa_{\mathcal{S}} (X) \to X$ has the right lifting property with respect to $I$ by the universal property of the arrow $\kappa_{\mathcal{S}}(X) \to X$. Hence it is a trivial fibration and so, in particular, a weak equivalence. It follows that the Quillen adjunction $i: \mathcal{M}_{\mathcal{S}} \leftrightarrows \mathcal{M}: \kappa_{\mathcal{S}}$ is a Quillen equivalence.     
 \end{proof}

\begin{remark} Note that property (1) in the definition of cofibrantly generated model categories from Section 2 is not needed in the last proof. This observation is due to Ji\v{r}\'{\i} Rosick\'{y}.
\end{remark}

\parindent=0.2in The assumption of Vop\v{e}nka's principle can be dropped as long as the local presentability of $\mathcal{M}_{\mathcal{S}}$ can be asserted otherwise. For example, if the objects in $S$ are presentable, then $\mathcal{M}_{\mathcal{S}}$ is locally presentable by Proposition \ref{proposition}. For many purposes, the condition that the objects in $S$ are presentable seems to be almost as good as knowing that $\mathcal{M}$ is a combinatorial model category.  An example of this will be shown in Proposition \ref{proposition2} below. \

\parindent=0.2in Moreover, $\mathcal{M}_{\mathcal{S}}$ is always locally presentable when $\mathcal{M}$ is a fibre-small topological category by [\cite{FR}, Theorem 3.6]. The category of topological spaces clearly has this property for example. 

\parindent=0.2in The following remarkable result, due to J. Smith,  is very useful for generating model category structures on locally presentable categories.

\begin{theorem} Let $\mathcal{C}$ be a locally presentable category, $I$ be a set of morphisms and $\mathcal{W}$ a class of morphisms. Suppose that the following are satisfied:
\begin{itemize}
\item[(i)] $\mathcal{W}$ satisfies the 2-out-of-3 property and it is closed under retracts in $\mathcal{C}^{\rightarrow}$,
\item[(ii)] $I-injective \subseteq \mathcal{W}$,
\item[(iii)] $Cof(I) \cap \mathcal{W}$ is closed under pushouts and transfinite compositions,
\item[(iv)] $\mathcal{W}$ is accessible and accessibly embedded in $\mathcal{C}^{\rightarrow}$. 
\end{itemize}
Then the classes $\mathcal{W}$, $Cof(I)$ and $ (Cof(I) \cap \mathcal{W})-injective$ define a combinatorial model category structure on $\mathcal{C}$.
\end{theorem}
\begin{proof}
This is a version of J. Smith's theorem [\cite{Be}, Theorem 1.7]. It reduces to it from the fact that an accessible and accessibly embedded subcategory of a locally presentable category is cone-reflective (see [\cite{AR}, Theorem 2.53]) and therefore it satisfies the solution-set condition at every object. 
\end{proof}

\parindent=0.2in Let $\mathcal{C}$ be a locally presentable category and $\mathcal{W}$ a class of weak equivalences that satisfies (i) and (iv) of the Theorem. The class of formal cofibrations $Cof$ with respect to the pair $(\mathcal{C}, \mathcal{W})$ is the largest cofibrantly closed class of morphisms such that $Cof \cap \mathcal{W}$ is closed under pushouts and transfinite compositions. Clearly, every cofibration of a model category
structure on $\mathcal{C}$ with weak equivalences $\mathcal{W}$ is a formal cofibration. Conversely, given a set $I$ of formal cofibrations that is big enough in the sense that $I-injective \subseteq \mathcal{W}$, then $Cof(I)$ is a class of cofibrations for such a (combinatorial) model structure. \

\begin{proposition} \label{proposition2}
\begin{itemize}
\item [(i)] Let $\mathcal{M}$ be a cofibrantly generated model category with $I$ and $J$ sets of generating cofibrations and trivial cofibrations respectively. Suppose that the domains and codomains of the maps in $I \cup J$ are presentable objects. Then the class of weak equivalences $\mathcal{W}$ is accessibly embedded in $\mathcal{M}^{\rightarrow}$.
\item [(ii)]Let $\mathcal{M}$ be a cofibrantly generated model category with $I$ and $J$ sets of generating cofibrations and trivial cofibrations respectively. Suppose that the codomains of the maps in $I$ are presentable objects.  Assuming Vop\v{e}nka's principle, the class of weak equivalences $\mathcal{W}$ is accessibly embedded in $\mathcal{M}^{\rightarrow}$. 
\end{itemize}
\end{proposition}
\begin{proof}
\begin{itemize}
\item [(i)] By Proposition \ref{proposition} (ii) and Theorem \ref{A}, there is a set of objects $S$ in $\mathcal{M}$ and a Quillen equivalence $\mathcal{M}_{\mathcal{S}} \leftrightarrows \mathcal{M}$ where $\mathcal{M}_{\mathcal{S}}$ is a cominatorial model category. By [\cite{D2}, Proposition 7.3]  , the class of weak equivalences $\mathcal{W}_{\mathcal{S}}$ of $\mathcal{M}_{\mathcal{S}}$ is accessibly embedded in $\mathcal{M}_{\mathcal{S}}^{\rightarrow}$, so they are closed under $\mu$-directed colimits for some $\mu$. Let $\lambda$ be a regular ordinal larger than $\mu$ such that every object that appears as the domain or codomain of $I$ is $\lambda$-presentable. We show that the weak equivalences $\mathcal{W}$ of $\mathcal{M}$ are closed under $\lambda$-directed colimits. Let $F: J \to \mathcal{M}^{\rightarrow}$ be a $\lambda$-directed diagram such that $F(j) \in \mathcal{W}$ for all $j$. Then $F_{\mathcal{S}}:= \kappa_{\mathcal{S}} \circ F: J \to \mathcal{M}_{\mathcal{S}}^{\rightarrow}$ satisfies $F_{\mathcal{S}}(j) \in \mathcal{W}_{\mathcal{S}}$ and therefore $colim_J F_{\mathcal{S}} \in \mathcal{W}_{\mathcal{S}}$. The arrows $F_{\mathcal{S}}(j) \to F(j)$ in $\mathcal{M}^{\rightarrow}$ are objectwise trivial fibrations by the proof of Theorem \ref{A}. Then the induced arrow $colim_J F_{\mathcal{S}} \to colim_J F$ in $\mathcal{M}^{\rightarrow}$ is also a trivial fibration objectwise because every member of $I$ is $\lambda$-presentable in $\mathcal{M}^{\rightarrow}$. It follows, by the 2-out-of-3 property, that $colim_J F \in \mathcal{W}$. Thus $\mathcal{W}$ is accessibly embedded in $\mathcal{M}^{\rightarrow}$.
\item [(ii)] The same argument essentially as in (i) applies to show that $\mathcal{W}$ is accessibly embedded in $\mathcal{M}^{\rightarrow}$. The assertion that $\mathcal{M}_{\mathcal{S}}$ is locally presentable requires Vop\v{e}nka's principle in this case. Also, the members of $I$ are not necessarily presentable in $\mathcal{M}^{\rightarrow}$ in this case, but they are presentable in the full subcategory of arrows in $\mathcal{M}^{\rightarrow}$ whose domains are in $\mathcal{M}_{\mathcal{S}}$, which suffices for the argument.  
\end{itemize}
\end{proof}

\section{The proof of Theorem \ref{B}}

\begin{proof} (of Theorem \ref{B})
Let $\mathcal{M}$ be a cofibrantly generated model category with generating set $I= \{A_i \to X_i \}_i$ of cofibrations such that $A_i$ is cofibrant for all $i$.  Let $1$ denote the terminal object in $\mathcal{M}$. Suppose that $X$ is an object in $Ho \mathcal{M}$ which we can assume to be fibrant in $\mathcal{M}$ and write $e: X \to 1$ for the unique morphism. Suppose that $e_*: Ho \mathcal{M}(K,X) \to Ho \mathcal{M}(K, 1) = \{ \star \}$ is a bijection for every $K \in \{A_i, X_i \}_i$. Then every diagram 
\begin{displaymath}
 \xymatrix{
A_i \ar[r]^f \ar[d]^i & X \ar[d] \\
X_i \ar[r] & 1
}
\end{displaymath}
admits a lift $h$ up to homotopy.  Since $A_i$ is cofibrant by assumption, there is a homotopy $H: Cyl(A_i) \cup_{A_i} X_i \to  X$ that restricts to $f$ and $h$ respectively. Note that $Cyl: \mathcal{M} \to \mathcal{M}$ denotes a functorial choice of cylinder objects in $\mathcal{M}$. We will show that $f$ extends to $X_i$.  There is a (cofibration, trivial fibration)-factorisation $j=p \circ t: Cyl(A_i) \cup_{A_i} X_i \to C' X_i \to Cyl(X_i)$. $j$ is weak equivalence, therefore so is $t$. The diagram
\begin{displaymath}
 \xymatrix{
Cyl(A_i) \cup_{A_i} X_i \ar[r] \ar[d]^t & X \ar[d] \\
C'X_i \ar[r] \ar@{.>}[ru]^g & 1 
}
\end{displaymath}
admits a lift $g: C'X_i \to X$. Also let $s: X_i \to C'X_i$ be a lift of 
\begin{displaymath}
 \xymatrix{
A_i \ar[r]^{i_1} \ar[d]^i & C'X_i \ar[d]^p \\
X_i \ar[r]^{i_i} \ar@{.>}[ru]^s & Cyl(X_i)
}
\end{displaymath}
Then $g \circ s: X_i \to X$ makes the first diagram commute. It follows that $e$ is a weak equivalence and $\{ A_i, X_i \}$ a left weakly adequate set of $Ho \mathcal{M}$. 
\end{proof}

\begin{remark} Note that the proof only requires that the class of cofibrations is cofibrantly generated by a set of cofibrations between cofibrant objects.
\end{remark}

\parindent=0.2in Str\o{}m's model category \cite{St} is a model category structure on the category $Top$ of topological spaces which is different from the usual cofibrantly generated model category of spaces \cite{Ho}. The cofibrations are the closed (Hurewicz) cofibrations (or closed NDR-pairs), the weak equivalences are the homotopy equivalences and the fibrations are the Hurewicz fibrations. With respect to this model category structure, every object is both cofibrant and fibrant. It follows that the homotopy category is the quotient of $Top$ with respect to the relation of homotopy on the morphism sets. It seems possible that it does not have a left weakly adequate set, but we do not know if this is true.
\

\parindent=0.2in Note that the stronger assertion that the homotopy category of a cofibrantly generated model category has a left adequate set is not true. For example, the usual cofibrantly generated model category of topological spaces does not have this property \cite{He}. \\ 

$\mathbf{Acknowledgements.}$ I would like to thank my supervisor Ulrike Tillmann for her support. I would also like to gratefully acknowledge the support of a partial EPSRC Studentship and a scholarship from the Onassis' Public Benefit 
Foundation.

\end{document}